\documentclass[12pt]{amsart}
\usepackage{graphicx} % Required for inserting images
\usepackage{amsmath,amsfonts,amssymb,amsthm,enumerate,color,cite,hyperref}

\def\R{\mathbb{R}}
\def\C{\mathbb{C}}
\def\A{\mathbb{A}}
\def\M{\mathcal{M}}
\def\dH{\dim_{\mathcal{H}}}

\newcommand{\supp}{\operatorname{supp}}

\newtheorem{thm}{Theorem}[section]

\newtheorem{definition}[thm]{Definition}

\renewcommand\Re{\operatorname{Re}}
\renewcommand\Im{\operatorname{Im}}

\title[]{Orthogonal projection, dual Furstenberg problem, and discretized sum-product}
\author[]{Longhui Li}
\address{Department of Mathematics \& International Center for Mathematics, Southern
University of Science and Technology, Shenzhen, 518055, PR China}
\email{12231267@mail.sustech.edu.cn}
\author[]{Bochen Liu}
\email{Bochen.Liu1989@gmail.com}

\date{}

\begin{document}

\begin{abstract}
    In this paper we come up with a dual version of the Furstenberg set problem and obtain partial results via $L^p$ estimates of orthogonal projections. Examples are also discussed. Moreover, compared with general sets, we find that special structure like Cartesian product has better $L^p$-behavior. This leads to improvement on some discretized sum-product estimates.
\end{abstract}

\maketitle

\section{Introduction}

For $x\in\R^2$ and $e\in S^1$, denote by $\pi_e(x):=x\cdot e$ the orthogonal projection. The classical Marstrand projection theorem \cite{Mar54} states that, for every Borel set $E$ in the plane of Hausdorff dimension $\dH E>1$, the set $\pi_e(E)\subset\R$ has positive Lebesgue
for almost every $e\in S^1$. Marstrand's original proof is very complicated. In 1968, Kaufman \cite{Kau68} gave a one-line proof using Fourier analysis:
\begin{equation}\label{Kaufman-one-line}\iint|\widehat{\pi_e\mu}(t)|^2\,dt\,de=\iint|\hat{\mu}(te)|^2\,dt\,de=c\int|\hat{\mu}(\xi)|^2\,|\xi|^{-1}d\xi.\end{equation}
To see this indeed implies the Marstrand projection theorem, one needs to take $\mu$ as a Frostman measure on $E$. More precisely, the well-known Frostman Lemma guarantees that, for every Borel subset $E\subset\R^d$ and every $\alpha<\dH E$, there exists a nonzero finite Borel measure $\mu$ supported on $E$ satisfying 
\begin{equation}\label{Frostman}
    \mu(B(x,r))\leq c r^\alpha, \ \forall\,x\in\R^d,\forall\,r>0.
\end{equation}
As a consequence, the $s$-dimensional energy, defined by
\begin{equation}\label{energy}I_s(\mu):=\iint |x-y|^{-s}\,d\mu(x)\,d\mu(y)=c_d\int_{\R^d}|\hat{\mu}(\xi)|^2\,|\xi|^{-d+s}\,d\xi,\end{equation}
is finite for every $s<\alpha<\dH E$. In addition to Marstrand's statement, Kaufman's proof also implies that for almost every $e\in S^1$ the induced measure $\pi_e\mu$ on $\pi_e(E)$ has $L^2$ density. Throughout this paper we call a measure $\mu$ satisfying \eqref{Frostman} a Frostman measure of exponent $\alpha$.
For a detailed discussion on these classical results including the Fourier-analytic representation of $I_s(\mu)$ in \eqref{energy}, we refer to \cite{Mat15}.

Notice Kaufman's one-line proof consists of equalities. This means nothing changes when $\dH E$ gets much larger than $1$. So one can ask about the benefit of taking $\dH E$ large. One point of view is, by Sobolev embedding and Kaufman's argument,
$$\int\|\pi_e\mu\|_{L^p(\R)}^2\,de\lesssim \int\|\pi_e\mu\|_{L^2_{(\frac{1}{2}-\frac{1}{p})}(\R)}^2\,de = c\int|\hat{\mu}(\xi)|^2\,|\xi|^{-2/p}d\xi.$$
As a consequence one can find $\mu$ on $E$ such that $\pi_e\mu\in L^p$ if $\dH E> 2(1-\frac{1}{p})$, or equivalently, $p<\frac{2}{2-\dH E}$. In \cite{PS00}, Peres and Schlag asked if this range of $p$ is optimal, but there is no further result. Notice it is a $(p,2)$ mixed-norm estimate on $\R\times S^1$.

There are analogs in higher dimensions, where the orthogonal projection is defined by $\pi_V:\R^d\mapsto V\approx \R^n$ for every $n$-dimensional subspace $V\in G(d,n)$. Here $G(d,n)$ denotes the Grassmannian with Haar measure $\gamma_{d,n}$. We again refer to \cite{Mat15} for details.

For a while there is no further discussion on $L^p$ estimates of orthogonal projections, possibly due to the lack of applications of the $(p,2)$ estimate. In \cite{DOV22}, Dabrowski, Orponen and Villa consider $(p,p)$ estimates and prove
\begin{equation}\label{Lp-DOV}\int\|\pi_V\mu\|^p_{L^p(\mathcal{H}^n)}\,d\gamma_{d,n}(V)<\infty,\ \forall\,p<2+\frac{s-n}{d-s},\end{equation}
given $\mu$ a Frostman measure of exponent $s>n$. Here the range of $p$ is in general sharp. Unlike Sobolev embedding, their estimate has applications in the Furstenberg set problem and discretized sum-product estimates. Later an alternative proof is given by the second author \cite{Liu24}, with a more explicit upper bound in terms of a newly defined quantity called $\alpha$-dimensional amplitude of $\mu$,
$$A_\alpha(\mu):=\sup_x \int |x-y|^{-\alpha}\,d\mu(y).$$ 
Both $s$-dimensional energy $I_s(\mu)$ and $\alpha$-dimensional amplitude $A_\alpha(\mu)$ give equivalent definitions of Hausdorff dimension, namely
$$\begin{aligned}\dH E=&\sup\{s: \exists\,\mu\in\mathcal{M}(E),\ I_s(\mu)<\infty\}\\=&\sup\{\alpha: \exists\,\mu\in\mathcal{M}(E),\ A_\alpha(\mu)<\infty\}.\end{aligned}$$
See \cite{Liu24} for detailed discussions. Here and throughout this paper $\M(E)$ denotes the space of finite Borel measures supported on a compact subset of $E$.

\subsection{$L^p$ estimates of orthogonal projections}
Denote by $\mathbb{A}(d,k)$ the affine Grassmannian consisting of all affine $k$-dimensional subspaces in $\R^d$. Notice that every element of $ \mathbb{A}(d,k)$ can be uniquely represented by $W+u$, where $W\in G(d,k)$ and $u\in W^\perp\subset\R^d$. Then there is a natural measure $\lambda_{d,k}$ on $\A(d,k)$ that is invariant under the group of rigid motions,
$$\begin{aligned}
\int_{\A(d,k)} f\,d\lambda_{d,k}:=&\int_{G(d,k)}\int_{W^\perp} f(W+u)\,d\mathcal{H}^{d-k}(u)\,d\gamma_{d,k}(W)\\=&\int_{G(d,d-k)}\int_{V} f(V^\perp+u)\,d\mathcal{H}^{d-k}(u)\,d\gamma_{d,d-k}(V).\end{aligned}$$
Here the last equality follows from the natural identification between $G(d,k)$ and $G(d, d-k)$, given by $V\mapsto V^\perp$.

With this measure $\lambda_{d,k}$ the estimate \eqref{Lp-DOV} can be written as
\begin{equation}\label{rewrite-Lp-DOV}\int_{\A(d,d-n)}|\pi_{V}\mu(u)|^p\,d\lambda_{d,d-n}(V^\perp+u)<\infty\end{equation}
because 
\begin{equation}\label{orthogonal-proj-integral-form}\pi_{V}\mu(u)=\int_{V^\perp+u}\mu\,d\mathcal{H}^{d-n}\end{equation} when $\mu$ has continuous density. In fact one can always assume $\mu\in C_0^\infty$, non-negative. We explain this in Section \ref{Sec-proof-Lp}.

In this paper we shall work with a more general measure defined by
\begin{equation}\label{measure-on-affine-subspaces}\begin{aligned}
    \int_{\mathbb{A}(d,k)} f\,d\lambda_{d,k,\nu}=&\int_{G(d,k)}\int_{W^\perp} f(W+u)\,d\mathcal{H}^{d-k}(u)\,d\nu(W)\\=&\int_{G(d,d-k)}\int_{V} f(V^\perp+u)\,d\mathcal{H}^{d-k}(u)\,d\nu(V),\end{aligned}
\end{equation}
where $\nu$ is a Frostman measure in $\M(G(d,k))=\M(G(d,d-k))$ of exponent $\sigma$. To study the orthogonal projection $\pi_V$ with $V\in G(d,n)$, one should take $k=d-n$ as in \eqref{rewrite-Lp-DOV}. Our $L^p$ estimates of orthogonal projections are the following. For convenience in the application we prefer writing $\int_{\A(d,d-n)}|\pi_V\mu(u)|^p\,d\lambda_{d,d-n,\nu}(V^\perp+u)$ instead of $\int_{G(d,d-n)}\|\pi_V\mu\|^p_{L^p}\,d\nu(V)$, though equivalent.
\begin{thm}\label{thm-Lp}
Suppose $\mu\in\M(\R^d)$ and $\lambda_{d,k,\nu}$ is defined as \eqref{measure-on-affine-subspaces} with $k=d-n$. Then, with $\alpha\in(0, d)$, $s>n(d-n+1)-\sigma$, and 
$$p:=2+\frac{s+\sigma-n(d-n+1)}{d-\alpha},$$ 
we have
$$\int_{\A(d,d-n)}|\pi_V\mu(u)|^p\,d\lambda_{d,d-n,\nu}(V^\perp+u)\lesssim I_s(\mu)\cdot A_\alpha(\mu)^{p-2}.$$
Moreover, if $\mu$ satisfies the Fubini property
\begin{equation}\label{mu-product-like}d\mu(x_1,x_2)=d\mu^{x_2}_1(x_1)\,d\mu_2(x_2),\ \mu^{x_2}_1\in\M(\R^n), \mu_2\in\M(\R^{d-n}),\end{equation}
and $span\{V, \{0\}\times\R^{d-n}\}=\R^d$ for every $V\in supp(\nu)\subset G(d,n)$, then with $0<\alpha<n$,  
$$p:=2+\frac{s+\sigma-n(d-n+1)}{n-\alpha},$$
we have
$$\int_{\A(d,d-n)}|\pi_V\mu(u)|^p\,d\lambda_{d,d-n,\nu}(V^\perp+u)\lesssim I_s(\mu)\cdot \sup_{x_2\in\supp\mu_2} A_\alpha(\mu^{x_2}_1)^{p-2}.$$
\end{thm}

An example satisfying the Fubini property \eqref{mu-product-like} is the Cartesian product $\mu=\mu_1\times\mu_2\in\M(\R^n)\times\M(\R^{d-n})$. In fact this result is motivated by considering Cartesian products. But it turns out we need this slightly more general version for the discretized sum-product estimates. See Theorem \ref{thm-A+B-AC} and Section \ref{sec-sum-product} below for details.

As a comparison, for general $\mu$ and particular $\nu=\lambda_{d,k}$ our range of $p$ coincides with that in \cite{DOV22}\cite{Liu24}, while the method here works on more general $\nu$ and obtains better results for special $\mu$. Also arguments in this paper are more straightforward.

%Another related example is to take $\mu=\lambda_G$ as a Haar measure on a locally compact group $G$ and $\mu_1^{x_2}=\lambda_H(x_2^{-1}\cdot)$ on a closed subgroup $H$ such that the modular function $\Delta$ satisfies $\Delta_G|_H=\Delta_H$. We refer to Theorem 2.51 in \cite{Fol16}. It would be very interesting if the idea in this paper works in a more general context.

\subsection{Incidence estimates}
Similar to \cite{DOV22}, for applications we first reduce our $L^p$ estimates to incidence estimates. Though this idea is not new, our proof is much shorter than \cite{DOV22}. 

More precisely, given Borel sets $E\subset\R^d$, $\mathcal{A}\subset\A(d,k)$ and $\delta>0$, we would like to figure out the size of
$$I_\delta(E,\mathcal{A}):=\{(x,W+u)\in E\times \mathcal{A}: x\in \mathcal{N}_\delta(W+u)\}.$$

\begin{thm}\label{thm-incidence}
    Suppose $E\subset\R^d$, $\mathcal{A}\subset\A(d,k)$ are Borel sets, $\mu\in\M(E)$ and $\lambda_{d,k,\nu}\in\M(\A(d,k))$ is defined as \eqref{measure-on-affine-subspaces}, then
    $$\mu\times\lambda_{d,k,\nu} (I_\delta(E,\mathcal{A}))\lesssim I_s(\mu)^\frac{1}{p}\cdot A_\alpha(\mu)^\frac{p-2}{p}\cdot\lambda_{d,k,\nu}(\mathcal{N}_{\delta}(\mathcal{A}))^\frac{1}{p'}\cdot\delta^{d-k}, $$
    where $$p:=2+\frac{s+\sigma-(k+1)(d-k)}{d-\alpha}. $$
    Moreover, if $\mu$ satisfies the Fubini property \eqref{mu-product-like} and $span\{V, \{0\}\times\R^k\}=\R^d$ for every $V\in supp(\nu)\subset G(d,d-k)$, then
    $$\mu_1\times\lambda_{d,k,\nu} (I_{\delta}(E,\mathcal{A}))\lesssim I_s(\mu)^\frac{1}{p}\cdot \sup_{x_2} A_\alpha(\mu^{x_2}_1)^{\frac{p-2}{p}}\cdot\lambda_{d,k,\nu}(\mathcal{N}_{\delta}(\mathcal{A}))^\frac{1}{p'}\cdot\delta^{d-k}, $$
    where $$p:=2+\frac{s+\sigma-(k+1)(d-k)}{d-k-\alpha}.$$
\end{thm}

Here we prefer writing the right hand side in terms of $I_s$ and $A_\alpha$, because in different applications the non-concentration conditions on $\mu$ may be different (e.g. \eqref{Frostman}, \eqref{non-concentration-condition-1}, \eqref{non-concentration-condition-2} in this paper).

\subsection{A dual version of the Furstenberg set problem}
From the incidence estimates in the last subsection, we can discuss a dual version of the Furstenberg set problem. For $s \in (0, 1]$ and $t\in (0, 2]$, an $(s, t)$-Furstenberg set in the plane is a set $E\subset\R^2$ such that there exists a family of lines $\mathcal{L}$ with $\dH \mathcal{L}\geq t$ satisfying $\dH(E\cap l)\geq s$ for all $l\in\mathcal{L}$. Accumulating efforts from many famous mathematicians, the minimal possible $\dH E$ has been fully understood recently. We refer to \cite{RW23} and references therein for the most recent work and the history of this problem. Their result is the following.
\begin{thm}
    An $(s, t)$-Furstenberg set in the plane has Hausdorff dimension at least $\min\{s+t, \frac{3s+t}{2}, s+1\}$.
\end{thm}
There are higher dimensional analogs that ask about $\dH E$ if there exists a set of $k$-planes whose intersection with $E$ has dimension at least $s$. For partial results we refer to \cite{DOV22,Hera19,HKM19}. 

In this paper we consider a dual version of the Furstenberg set problem, by exchanging the roles of points and $k$-planes. More precisely, given a subset $\mathcal{A}\subset\A(d,k)$, suppose there exists a set $E\subset\R^d$ with $\dH E\geq s$  such that for every $x\in E$,
$$\dH\{W+u\in \mathcal{A}: x\in W+u\}\geq t,$$
then what is the minimal possible $\dH \mathcal{A}$? To define Hausdorff dimension in the ambient space $\A(d,k)$ one can consider the usual metric (see, e.g. Section 3.16 in \cite{Mat95}) given by
$$d_{\mathbb{A}(d,k)}(W+u, W'+u'):=d_{G(d,k)}(W, W')+|u-u'|=\|\pi_W- \pi_{W'}\|+|u-u'|.$$
It seems people from other areas may prefer the metric defined by
$$d_{G(d+1,k+1)}(span\{(W,0), (u,1)\}, span\{(W',0), (u',1)\}),$$
where $W+u\mapsto span\{(W,0), (u,1)\}$ is the natural embedding from $\mathbb{A}(d,k)$ to $G(d+1, k+1)$ (see, e.g. \cite{LWY21} and references therein). But for our use these two metrics are equivalent, as one can always reduce the set $E$ above to a bounded subset $E'$ (see Section \ref{sec-dual-Furstenberg} below), and then it suffices to consider affine subspaces that intersect $E'$. Moreover, in this case $d(W+u, W'+u')<\delta$ if and only if $E'\cap (W+u)$ is contained in the $c\delta$-neighborhood of $W'+u'$. Also notice that 
$$\dim \A(d,k)=\dim G(d,k) + \dim \R^{d-k}=(k+1)(d-k),$$
which equals $d$ if and only if $k=d-1$. This means the case $k=d-1$ is very special. In fact, when $k=d-1$ our dual version is equivalent to the original Furstenberg set problem by applying the projective transformation, also known as the point-plane duality. But for $k<d-1$ it seems to be an independent problem that has its own interest. We call the set $\mathcal{A}$ above an $(s,t)$-Furstenberg set in $\A(d,k)$.
\begin{definition}
A subset $\mathcal{A}\subset\A(d,k)$ is called an $(s,t)$-Furstenberg set in $\A(d,k)$, if there exists a set $E\subset\R^d$ with $\dH E\geq s$  such that
$$\dH\{W+u\in \mathcal{A}: x\in W+u\}\geq t,\ \forall\,x\in E.$$   
\end{definition}
One can also add the assumption that $E$ or $\mathcal{A}$, or both, are Borel or even compact sets, but in this paper we shall run a very careful discretization argument to avoid any of these assumptions. 

To state our results in full generality, we also need the natural projection $p$ from $\A(d,k)$ to $G(d,k)$, i.e., $p(W+u):=W$.

\begin{thm}\label{thm-dual-Furstenberg}
Under notation above we have the following.
\begin{enumerate}[(1)]
    \item Suppose $\mathcal{A}$ is an $(s,t)$-Furstenberg set in $\A(d,k)$ and $\mathcal{N}_\delta(p(\mathcal{A}))$ is a $(\delta, \sigma)$-set in the sense of \eqref{non-concentration-condition-2}  below for every $\delta\in(0,1)$, $s>(k+1)(d-k)-\sigma$, then
    $$\dH \mathcal{A}\geq t+(d-k)-\frac{(d-s)(\sigma-t)}{d+\sigma-(k+1)(d-k)}.$$
    \item If, in addition, the set $E$ can be taken as a Cartesian product $E_1\times E_2\subset\R^{d-k}\times\R^k$, $0<\mathcal{H}^{s_1}(E_1), \mathcal{H}^{s_2}(E_2)<\infty, s_1+s_2=s$, then the dimensional exponent can be improved to
    $$\dH \mathcal{A}\geq t+(d-k)-\frac{(d-k-s_1)(\sigma-t)}{d-k+s_2+\sigma-(k+1)(d-k)}.$$
\end{enumerate}
\end{thm}

Here the assumption on $\pi(\mathcal{A})$ does not lose much generality: one can always take $\sigma=\dim G(d,k)=k(d-k)$, in which case the dimensional exponents become
$$\dH \mathcal{A}\geq 2t-(k-1)(d-k)+\frac{(s-(d-k))(k(d-k)-t)}{k}$$
and
$$\dH \mathcal{A}\geq 2t-(k-1)(d-k)+\frac{(s-(d-k))(k(d-k)-t)}{s_2}.$$
Also $\mathcal{N}_\delta(p(\mathcal{A}))$ is guaranteed to be $(\delta, \sigma)$ if $p(\mathcal{A})$ is Alhfors-David regular. We hope that the assumption on $p(\mathcal{A})$ can be relaxed to $\dH p(\mathcal{A})=\sigma$ in later work.

The assumption $0<\mathcal{H}^{s_1}(E_1), \mathcal{H}^{s_2}(E_2)<\infty$ in the second part of Theorem \ref{thm-dual-Furstenberg} is also not very strong. In fact every Borel set $E$ with $\dH E>s$ has a compact subset $E'$ with $0<\mathcal{H}^s(E')<\infty$ (see, e.g, Chapter 8 in \cite{Mat95}). We shall explain why this extra condition is needed on Cartesian products in Section \ref{subsec-dual-Furstenberg-product}.

Compared with previous work, when $k=d-1$ the dual version is equivalent to the original Furstenberg set problem. In this case our results generalize those in \cite{DOV22} with the parameter $\sigma$ and improve those in \cite{DOV22} on Cartesian products. When $k< d-1$ these are the first results on the dual Furstenberg set problem.

Now we discuss potential sharp examples. Suppose $E\subset \mathbb{R}^d$, $\dH E=s$ and $\mathcal{V}\subset G(d,n)$, $\dH\mathcal{V}=t$, then
$$\mathcal{A}:=\{V+\pi_V(x): V\in\mathcal{V}, x\in E\}$$
is clearly a $(s,t)$-Furstenberg set in $\mathcal{A}(d,n)$. This means one can construct examples on our dual Furstenberg set problem from orthogonal projections. For dimension of $\mathcal{A}$, one may expect 
\begin{equation}\label{expected-upper-bound}\dH \mathcal{A}\leq \dH\mathcal{V}+\sup_{V\in\mathcal{V}}\dH\pi_V E,\end{equation}
but it is not guaranteed in general. Fortunately all potential sharp examples we know on orthogonal projections (see, e.g. \cite{KM75}\cite{Gan24}) are constructed by taking intersection of $\delta$-neighborhoods of arithmetic progressions, with \eqref{expected-upper-bound} holds as a trivial dimensional upper bound. We believe these are the optimal. For general $k$ the precise formula seems complicated, while in particular when $d=2, k=1$ it becomes
\begin{equation}\label{lower-bound-k=d-1}\dH\mathcal{A}\leq\min\left\{t+s, \frac{3t+s}{2}, t+1\right\}.\end{equation}
By the point-plane duality these become sharpness examples of $(t,s)$-Furstenberg sets in the plane (see \cite{RW23}). Previously such sharpness examples are only known from Wolff's direct construction (see \cite{Wolff99Kakeya}). After comparison, it seems easier to obtain these dimensional exponents from the dual version, as it is more convenient to compute projections than intersections. Another related remark is Bourgain's proof on the planar Kakeya also relies on considering orthogonal projections \cite{Bou91}.

%As the end of this subsection we would like to explain why here we can only study the dual version but not the original Furstenberg (except the special case $k=d-1$). The main reason is, to make our incidence estimates explicit, we need non-concentration conditions on $\mu$ but not  necessarily on $\lambda$. So it can only be used to consider arbitrary coverings of $\mathcal{A}\subset\A(d,k)$ with a given Frostman measure $\mu$ on $E\subset\R^d$, not the other way around. In the special case $k=d-1$, the roles of points and planes are exchangable, so there are more aspects to study this problem. As an example, Bourgain proved the planar Kakeya by considering orthogonal projections. %see e,g, \cite{Wol03}.
%This is crucial to us when computing $\dH(\mathcal{A})$ because by the definition of $\dH$ we need to work with an arbitrary covering of $\mathcal{A}$. For this reason our incidence estimates help only on the dual version of Furstenberg set problem to compute $\dH\mathcal{A}$, but not on the original version to compute $\dH E$. This also explains that why sharp incidence estimates on $\alpha$-dimsntional balls and $\beta$-dimensional tubes, proved earlier by FR, could not fully solve the Furstenberg set problem in the plane. It is similar in higher dimensions.

\subsection{Discretized sum-product}
As we have seen above, Cartesian product, or more generally the Fubini property \eqref{mu-product-like}, has extra benefits. We shall take use of it in the study of discretized sum-product.

Let $A\subset[1,2]$ be a union of $\delta$-intervals. Under some non-concentration condition, one should expect that $\max\{|A+A|, |AA|\}$ is much larger than $|A|$. There are two common ways to describe the non-concentration: 
\begin{enumerate}[(i)]
    \item for every $\delta\leq r\leq 1$,
    \begin{equation}\label{non-concentration-condition-1}
        |A\cap [x-r, x+r]|\lesssim r^{s_A}|A|;
    \end{equation}
    \item for every $\delta\leq r\leq 1$,
    \begin{equation}\label{non-concentration-condition-2}|A\cap [x-r, x+r]|\lesssim \delta\cdot (r/\delta)^{s_A}.\end{equation}
\end{enumerate}
Notice that
\begin{itemize}
    \item if $A$ satisfies \eqref{non-concentration-condition-1}, then $\mu=\frac{1}{|A|}\chi_A$ is a Frostman measure and for $s,\alpha<s_A$,
$$I_s(\mu)\lesssim 1,\quad A_\alpha(\mu)\lesssim 1;$$
\item if $A$ satisfies \eqref{non-concentration-condition-2}, then for $s,\alpha<s_A$ one can easily check that
$$I_s(\chi_A)\lesssim \delta^{1-s_A}|A|,\quad A_\alpha(\chi_A)\lesssim \delta^{1-s_A}.$$
\end{itemize}
For convenience in this paper we assume $|A|\approx \delta^{1-s_A}$, in which case these two conditions are the same. We call such $A$ a $(\delta,s_A)$-set.
Of course one can only assume one of \eqref{non-concentration-condition-1}\eqref{non-concentration-condition-2} and obtain corresponding estimates.

The discretized sum-product problem was first studied by Katz-Tao \cite{KT01} and Bourgain \cite{Bou03}. The best known result is due to Ren-Wang \cite{RW23} and Fu-Ren \cite{FR24}, but still not fully solved yet. 

There are also other discretized sum-product estimates. For example one can consider $\max\{|A+B|, |AC|\}$, where $A, B, C$ are unions of $\delta$-intervals and some (not necessarily all) satisfy non-concentration conditions. It is already shown in \cite{DOV22} that $L^p$ estimates of projections imply partial results on $\max\{|A+B|, |AC|\}$. We would like to point it out that in this application projection estimates are used on sets with product-like structure, while no difference between general sets and Cartesian products was shown in \cite{DOV22}. This is where we obtain improvement.

\begin{thm}\label{thm-A+B-AC}
Suppose $A\subset[1,2]$ is a union of $\delta$-intervals, $B\subset[1,2]$ is a $(\delta, s_B)$-set, $C\subset[1,2]$ is a $(\delta, s_C)$-set, $s_B+s_C>1$, then for any $\epsilon>0$,
$$\max\{|A+B|, |AC|\}\gtrsim_{\epsilon, s_B, s_C} |A|^{1-\frac{s_B+s_C-1}{2\min\{s_B, s_C\}}+\epsilon}.$$
    In particular, by taking $B=C$, 
$$\max\{|A+B|, |AB|\}\gtrsim_{\epsilon, s_B} |A|^{\frac{1}{2s_B}+\epsilon}.$$
\end{thm}
For general $A, B, C$ this result is always better than that in \cite{DOV22} where the exponent is $1-\frac{s_B+s_C-1}{2}+\epsilon$. When $|A|=\delta^{1-s_A}$ and $B=C$ we can conclude that
$$\max\{|A+B|, |AB|\}\gg |B|$$
if $s_A>1-2s_B(1-s_B)$, which is also better than that in \cite{DOV22}. Here $\gg |B|$ means $\gtrsim \delta^{-\beta}|B|$ for some $\beta=\beta(s_A, s_B)>0$. Notice that our method, as well as that in \cite{DOV22}, gives new results on sum-product only if $A,B,C$ are quite different. When $A=B=C$ the best known result is still due to \cite{FR24}\cite{RW23}. 

\subsection*{Organization.}
This paper is organized as follows. In Section \ref{Sec-proof-Lp} we prove our $L^p$ estimates of orthogonal projections (Theorem \ref{thm-Lp}). In Section \ref{sec-incidence} we use projection estimates to prove incidence estimates (Theorem \ref{thm-incidence}). In Section \ref{sec-dual-Furstenberg} we use incidence estimates to obtain dimensional lower bound on dual Furstenberg set problem (Theorem \ref{thm-dual-Furstenberg}). In Section \ref{sec-sum-product} we use incidence estimates on measures with Fubini property to obtain new results on discretized sum-product (Theorem \ref{thm-A+B-AC}).
\subsection*{Notation.}
Throughout this article we write $A\lesssim B$ if there exists an absolute constant $C$ such that $A\leq CB$.

\section{$L^p$ estimates of orthogonal projections}\label{Sec-proof-Lp}
In this section we prove Theorem \ref{thm-Lp}. Suppose $\phi\in C_0^\infty$, $\phi\geq 0$, $\int \phi =1$, $\phi_\epsilon(\cdot):=\epsilon^{-d}\phi(\frac{\cdot}{\epsilon})$, and let
$$\mu^\epsilon(x):=\phi_\epsilon*\mu(x).$$
Then the Riesz potential,
\begin{equation}
	\label{Riesz-potential}
	\mu^\epsilon_z(x):=\frac{\pi^\frac{z}{2}}{\Gamma(\frac{z}{2})}\,|\cdot|^{-d+z}*\mu^\epsilon(x),
\end{equation}
is a smooth function in $x\in\R^d$, initially defined for $\Re z>0$ and can be extended to $z\in\C$ by analytic continuation. It is well known that \eqref{Riesz-potential} has Fourier transform
$$\frac{\pi^\frac{d-z}{2}}{\Gamma(\frac{d-z}{2})}\,\hat{\mu^\epsilon}(\xi)|\xi|^{-z}$$
in the sense of distribution and in particular \eqref{Riesz-potential} equals $c_d\,\mu^\epsilon$ when $z=0$ as a distribution. We refer to \cite{GS77}, p71, p192, for details. An observation is, when $0<\Re z<d$,
\begin{equation}\label{ordinary-meaning-amplitude}\left\||\cdot|^{-d+z}*\mu^\epsilon\right\|_{L^\infty}\leq \left\||\cdot|^{-d+\Re z}*\mu^\epsilon\right\|_{L^\infty}=A_{\Re z}(\mu^\epsilon).\end{equation}

It is routine to consider $\mu^\epsilon$ first and take $\epsilon\rightarrow 0$ at the very end. One can check that all implicit constants below are independent in $\epsilon$. 
Also it is not hard to verify that
$$\lim_{\epsilon\rightarrow 0}A_\alpha(\mu^\epsilon)= A_\alpha(\mu),\ \lim_{\epsilon\rightarrow 0}I_s(\mu^\epsilon)=I_s(\mu),$$
if finite. From now we write $\mu$ for $\mu^\epsilon$ for abbreviation and assume $\mu\in C_0^\infty$. In this case the projection $\pi_V\mu(u)$ can be written as
\begin{equation}\label{proj-continuous-function}\pi_V\mu(u)=\int_{\pi_V^{-1}(u)}\mu\,d\mathcal{H}^{d-n}=\int_{\pi_V^{-1}(u)}\mu\,\psi\,d\mathcal{H}^{d-n}\end{equation}
for some $\psi\in C_0^\infty$ which is non-negative satisfying $\psi=1$ on $\supp\mu$. Fix this $\psi$ in this section.

Now we start our proof of Theorem \ref{thm-Lp}. For the first part, we are working with 
$$p:=2+\frac{s+\sigma-n(d-n+1)}{d-\alpha}=\frac{2(d-\alpha)+s+\sigma-n(d-n+1)}{d-\alpha}.$$
whose conjugate is
$$p'=1+\frac{d-\alpha}{(d-\alpha)+s+\sigma-n(d-n+1)}=\frac{2(d-\alpha)+s+\sigma-n(d-n+1)}{(d-\alpha)+s+\sigma-n(d-n+1)}.$$

For an arbitrary $f\in L^{p'}(\lambda)$ with $\|f\|_{L^{p'}}=1$, consider the analytic function
$$\Phi(z):=e^{z^2}\int_{G(d,n)}\int_V \pi_V(\psi\mu_z)(u)\cdot f_z(V^\perp+u)\,d\mathcal{H}^n(u)\,d\nu(V),$$
where $\mu_z$ is defined as in \eqref{Riesz-potential} and
$$f_z(V,y):= sgn(f) \cdot |f(V,y)|^{p'\frac{(d-\alpha-z)+(2z+s+\sigma-n(d-n+1))}{2(d-\alpha)+s+\sigma-n(d-n+1)}}.$$
The exponent over $|f|$ looks complicated but the point is to ensure it is linear in $z$ and equals
$$\begin{cases}
    p',& \text{when }z=d-\alpha>0;\\
    1,& \text{when }z=0;\\
    \frac{p'}{2},& \text{when }z=-\frac{s+\sigma-n(d-n+1)}{2}<0.
\end{cases}$$
Also the role of $e^{z^2}$ here is to ensure
\begin{equation}\label{control-Gamma-function}|e^{z^2}\cdot\Gamma^{-1}(z)|\lesssim_{Re z} 1,\end{equation}
as $\Gamma^{-1}$ is an entire function of order $1$.

Then it suffices to prove
$$|\Phi(0)|\lesssim I_s(\mu)^{1/p}\cdot A_\alpha(\mu)^{1-\frac{2}{p}}.$$

When $\Re z=d-\alpha>0$, we have $\|f_z\|_{L^1}\leq \|f_{\Re z}\|_{L^{p'}}=1$ and therefore
$$|\Phi(z)|\lesssim\|e^{z^2}\pi_V(\psi\mu_z)\|_{L^\infty}\lesssim_\psi\|e^{z^2}\mu_z\|_{L^\infty}\lesssim A_\alpha(\mu),$$
where all implicit constants are uniform in $\Im z$. Here the second inequality follows from the expression \eqref{proj-continuous-function} and the last inequality follows from \eqref{ordinary-meaning-amplitude}, \eqref{control-Gamma-function}.

When $\Re z=-\frac{s+\sigma-n(d-n+1)}{2}<0$, we have $\|f_z\|_{L^2}\leq \|f_{\Re z}\|_{L^{p'}}=1$ and therefore by Cauchy-Schwarz
$$|\Phi(z)|\lesssim\|\pi_V(e^{z^2}\psi\mu_z)\|_{L^2}.$$

Now we invoke a classical estimate that dates back to \cite{Fal82}, which is, under our notation,
$$\int|\pi_V\mu(u)|^2\,d\lambda(V^\perp+u)\lesssim I_{n(d-n+1)-\sigma}(\mu).$$
We also refer to \cite{Mat15}, Section 5.3 for detailed discussion. We would like to point out that, from the proof it is very easy to see that it remains valid when $\mu$ is complex-valued. Therefore it applies to $e^{z^2}\psi\mu_z$ and
$$|\Phi(z)|\lesssim\|\pi_V(e^{z^2}\psi\mu_z)\|_{L^2}\lesssim I_{n(d-n+1)-\sigma}(e^{z^2}\psi\mu_z)^{1/2},$$
which is
$$\lesssim I_{n(d-n+1)-\sigma-2\Re z}(\mu)^{1/2}=I_s(\mu)^{1/2}$$
by standard argument (see, e.g., Subsection 4.3 in \cite{Liu24} for details). Again all implicit constants here are uniform in $\Im z$.

Finally, we apply Hadamard's three-lines lemma to the analytic function
$$\frac{\Phi(z)}{I_s(\mu)^\frac{d-\alpha-z}{2(d-\alpha)+s+\sigma-n(d-n+1)}\cdot A_\alpha(\mu)^\frac{2z+s+\sigma-n(d-n+1)}{2(d-\alpha)+s+\sigma-n(d-n+1)}}$$
in the strip $$\left\{-\frac{s+\sigma-n(d-n+1)}{2}<\Re z<d-\alpha\right\}$$ to conclude that
$$|\Phi(0)|\lesssim I_s(\mu)^\frac{d-\alpha}{2(d-\alpha)+s+\sigma-n(d-n+1)}\cdot A_\alpha(\mu)^\frac{s+\sigma-n(d-n+1)}{2(d-\alpha)+s+\sigma-n(d-n+1)}=I_s(\mu)^{1/p}\cdot A_\alpha(\mu)^{1-\frac{2}{p}},$$
as desired.

The proof of the second part of Theorem \ref{thm-Lp} is quite similar. The only difference is, instead of $\mu_z$, we consider $(\mu_1)_z\times\mu_2$ and set
$$\Phi(z):=e^{z^2}\int_{G(d,n)}\int_V \pi_V((\mu^{x_2}_1)_z(x_1)\mu_2(x_2))(u)\cdot f_z(V^\perp+u)\,d\mathcal{H}^n(u)\,d\nu(V).$$
Thanks to the Fubini property and the assumption $span\{V, \{0\}\times\R^{d-n}\}=\R^d$, we have that, when $\Re z=n-\alpha$,
$$\begin{aligned}\|e^{z^2}\pi_V((\mu^{x_2}_1)_z(x_1)\mu_2(x_2)\|_{L^\infty}=&\left\|e^{z^2}\int_{(x_1,x_2)\in V^\perp+u}(\mu^{x_2}_1)_z(x_1)\,\mu_2(x_2)\,d\mathcal{H}^{d-n}\right\|_{L^\infty}\\\lesssim &\sup_{x_2}A_\alpha(\mu_1^{x_2})\cdot \left\|\int_{(x_1,x_2)\in V^\perp+u}\mu_2(x_2)\,d\mathcal{H}^{d-n}\right\|_{L^\infty}\\\lesssim &\sup_{x_2}A_\alpha(\mu_1^{x_2})\cdot\|\mu_2\|_{L^1(\R^{d-n})}.\end{aligned}$$
Then the rest of the proof is the same as the first part with $d-\alpha$ replaced by $n-\alpha$. We omit details.

\section{From $L^p$ estimates to incidence estimates}\label{sec-incidence}
In this section we prove Theorem \ref{thm-incidence}. First notice
$$(\mu\times\lambda_{d,k,\nu})(I_\delta(E,\mathcal{A}))=\int_{\mathcal{A}}\int_{x\in \mathcal{N}_\delta(W+u)}\mu(x)\,dx\,d\lambda_{d,k,\nu}(W+u).$$
For each fixed $W$ we change variables in $x\in\R^d$ in terms of $\R^d=W\oplus W^\perp$. Then by Fubini and the integral form of orthogonal projections \eqref{orthogonal-proj-integral-form}, this integral equals
$$\int_{\mathcal{A}}\int_{|u'-u|\leq\delta}\pi_{W^\perp}\mu(u')\,du'\,d\lambda_{d,k,\nu}(W+u)$$$$ =\iint_{W+u\in\mathcal{A}}\int_{|u'-u|\leq\delta}\pi_{W^\perp}\mu(u')\,du'\,du\,d\nu(W).$$
By H\"older's inequality, it is bounded from above by
    $$\begin{aligned}& \left(\iint_{|u-u'|\leq\delta}|\pi_{W^{\bot}}\mu(u')|^p\,du\,du'\,d\nu(W)\right)^{1/p}\cdot\left(\iiint_{\substack{u+W\in\mathcal{A},\\ |u-u'|\leq\delta}}du'\,du\,d\nu(W)\right)^{1/p'}\\\lesssim &\delta^{d-k}\cdot \left(\iint|\pi_{W^{\bot}}\mu(u')|^p\,du'\,d\nu(W)\right)^{1/p}\cdot\left(\iint_{u+W\in\mathcal{N}_\delta(\mathcal{A})}\,du\,d\nu(W)\right)^{1/p'}\\=&\delta^{d-k}\cdot \left(\iint|\pi_{W^{\bot}}\mu(u')|^p\,du'\,d\nu(W)\right)^{1/p}\cdot\lambda_{d,k,\nu}(\mathcal{N}_\delta(\mathcal{A}))^{1/p'}.\end{aligned}$$
Finally we apply the first part of Theorem 1.1 with $k=d-n$ to conclude the first part of Theorem 1.2. The second part of Theorem 1.2 follows in a similar way.

\section{From incidence to dual Furstenberg}\label{sec-dual-Furstenberg}
In this section we prove Theorem \ref{thm-dual-Furstenberg}. As none of our sets is assumed Borel,  we need a very careful discretization argument before employing our incidence estimates. We shall discuss two parts of the theorem separately as the product case is trickier.

\subsection{General sets}

For the first part, since $\dH E\geq s$, for every $0<s'<s$ there exists a bounded subset $E'\subset E$ with $\dH E'>s'$. Then its $s'$-dimensional Hausdorff content satisfies $$0<\mathcal{H}^{s'}_\infty(E')<\infty. $$
    Let $t'<t$ be arbitrary. For every $x\in E'$, by our assumption
  $$\mathcal{H}^{t'}_{\infty}(\{W+u\in\mathcal{A}:x\in W+u\})>0.$$
  By the sub-additivity of the Hausdorff content, there exists a constant $c>0$, a subset $E''\subset E'$ of $\mathcal{H}^{s'}_{\infty}(E'')=c'\in(0,\infty)$, such that
  $$\mathcal{H}^{t'}_{\infty}(\{W+u\in\mathcal{A}:x\in W+u\})>c>0,\ \forall\,x\in E''.$$
  From now $c>0$ and $c'\in (0,\infty)$ are fixed constants.
    
 Now suppose $\mathcal{H}^\gamma(\mathcal{A})=0$. This means for every large enough $j_0\in\mathbb{Z}_+$, there exists a countable cover $\left\{B_i\right\}_{i\in I}$  of $\mathcal{A}$ by balls, with $\mathrm{diam}\; B_i<2^{-j_0}$ for each $i\in I$, such that
 $$\sum (\mathrm{diam}\; B_i)^\gamma\leq 1.$$
 	Denote 
  $$\mathcal{B}_j:=\left\{B_i:2^{-j-1}\leq\mathrm{diam}\; B_i<2^{-j}\right\},\ j\geq j_0.$$

  By pigeonholing and the sub-additivity of the Hausdorff content, one can conclude that for every $x\in E''$ there exists $j\geq j_0$ such that 
  $$\mathcal{H}^{t'}_{\infty}(\{W+u\in\mathcal{A}\;\cap(\bigcup_{B_i\in\mathcal{B}_j}B_i):x\in W+u\})\geq \frac{\pi^2c}{6}\cdot \frac{1}{j^2}.$$ 
  As $\mathcal{B}_j$ is already a cover of $\mathcal{A}\;\cap(\bigcup_{B_i\in\mathcal{B}_j}B_i)$, this implies
  \begin{equation}\label{lower-bound-B-i}\#(\mathcal{B}_{j,x}):=\#\{B_i\in\mathcal{B}_j: \exists\,W+u\in B_i \text{ s.t. } x\in W+u\}\geq \frac{\pi^2c}{6\cdot 2^{t'}}\cdot \frac{2^{t'j}}{j^2}\end{equation}
  that will be useful later.

  Next we let $j_x$ denote the smallest such $j$ for each $x\in E''$. Then $$E''=\bigcup_{j\geq j_0}\left\{x\in E'':j_x=j\right\}$$
  and by pigeonholing there exists $j_1\geq j_0$ such that 
	$$\mathcal{H}^{s'}_\infty(\{x\in E'':j_x=j_1\})\geq \frac{\pi^2c'}{6}\cdot \frac{1}{j_1^2}.$$

 Now we denote $\delta=2^{-j_1}$,
 $$E_\delta:=\mathcal{N}_{\delta}(\{x\in E'':j_x=j_1\}),\  \mathcal{A}_\delta:=\bigcup_{B_i\in\mathcal{B}_{j_1}}B_i$$
 and consider their incidence. The reason we take the neighborhood is $\{x\in E'':j_x=j\}$ itself is not guaranteed to be a Borel set. Then, since $E_\delta$ is a bounded open set and $$\infty>\mathcal{H}^{s'}_{\infty}(E_\delta)\geq \mathcal{H}^{s'}_\infty(\{x\in E'':j_x=j_1\})\geq \frac{\pi^2c'}{6}\cdot \frac{1}{j_1^2},$$
 by the famous Frostman Lemma (see, e.g. Theorem 8.8 in \cite{Mat95}) one can find a finite Borel measure $\mu\in\M(E_\delta)$ such that
 $$\mu(B(x,r))\lesssim \mathcal{H}^{s'}_{\infty}(E_\delta)^{-1}\cdot r^{s'}\lesssim j_1^2\cdot r^{s'},$$
 where the implicit constant only depends on fixed constants $c,c',d$.

 Then we estimate 
 $$\mu\times\lambda_{d,k}(I_{3\delta}(E_\delta,\mathcal{A}_\delta)).$$

For lower bound, by our construction every $x'\in E_\delta$ is $\delta$-close to some $x$ with $j_x=j_1$, which means
\eqref{lower-bound-B-i} holds for $x$ with $j=j_1$. From this fact one can conclude that such $x'$ is $(3\delta)$-close to every element in every $B\in\mathcal{B}_{j_1, x}$. Therefore
\begin{equation}\label{lower bound of incidence}\mu\times\lambda_{d,k}(I_{3\delta}(E_\delta,\mathcal{A}_\delta))\gtrsim j_1^{-2}\cdot \delta^{-t'+(k+1)(d-k)}.\end{equation}
		
Now we turn to the upper bound. Recall we assume every $\delta$-neighborhood of $p(\mathcal{A})$, the direction set of $\mathcal{A}$, is $(\delta,\sigma)$ in the sense of \eqref{non-concentration-condition-2}. Therefore
$$\nu:=\delta^{-k(d-k)+\sigma}\cdot\gamma_{d,k}|_{\pi(\mathcal{A}_\delta)}$$ is a Frostaman measure of exponent $\sigma$. Consequently we have
	\begin{align*}
		\mu\times\lambda_{d,k}(I_{3\delta}(E_\delta,\mathcal{A}_\delta))\lesssim\delta^{k(d-k)-\sigma}\cdot\mu\times\lambda_{d,k,\nu}(I_{3\delta}(E_\delta,\mathcal{A}_\delta)).
	\end{align*}
	Then we are ready to invoke the first part of Theorem 1.2, with $$p=2+\frac{s'+\sigma-(k+1)(d-k)}{d-s'}\left(\implies p'=1+\frac{d-s'}{d+\sigma-(k+1)(d-k)}\right),$$ to get the upper bound
	\begin{equation}\label{upper bound of incidence}
 \begin{aligned}
		\mu\times\lambda_{d,k}(I_{3\delta}(E_\delta,\mathcal{A}_\delta))\lesssim & \delta^{(k+1)(d-k)-\sigma}\cdot I_{s'}(\mu)^{\frac{1}{p}}\cdot A_{s'}(\mu)^{\frac{p}{p-2}}\cdot \lambda_{d,k,\nu}(\mathcal{N}_{3\delta}(\mathcal{A}_\delta))^{\frac{1}{p'}}\\\lesssim & \delta^{(k+1)(d-k)-\sigma}\cdot \delta^{\frac{-k(d-k)+\sigma}{p'}}\cdot\lambda_{d,k}(\mathcal{N}_{3\delta}(\mathcal{A}_\delta))^{\frac{1}{p'}}\\\lesssim & \delta^{(k+1)(d-k)-\sigma}\cdot \delta^{\frac{d-k+\sigma}{p'}}\cdot\#(\mathcal{B}_{j_1})^{\frac{1}{p'}}
	\end{aligned}
 \end{equation}
	Combine (\ref{lower bound of incidence}) and (\ref{upper bound of incidence}), we have
 \begin{equation}\label{count-B-j-1}\#(\mathcal{B}_{j_1})\gtrsim j_1^{-2p'}\cdot\delta^{(\sigma-t')p'-(d-k+\sigma)}=j_1^{-2p'}\cdot\delta^{-t'-(d-k)+\frac{(d-s')(\sigma-t')}{d+\sigma-(k+1)(d-k)}}.\end{equation}

 Recall at the very beginning of the proof we assume $\mathcal{H}^\gamma(\mathcal{A})=0$ and
 $$\sum (\mathrm{diam}\; B_i)^\gamma\leq 1.$$
 So \eqref{count-B-j-1} implies
$$(-\log\delta)^{-2p'}\cdot\delta^{\gamma-t'-(d-k)+\frac{(d-s')(\sigma-t')}{d+\sigma-(k+1)(d-k)}}\lesssim 1,$$
where $\delta\leq 2^{-j_0}$ and the implicit constant is independent in $j_0$. Hence we can take $j_0\rightarrow\infty$ to obtain 
$$\gamma\geq t'+(d-k)-\frac{(d-s')(\sigma-t')}{d+\sigma-(k+1)(d-k)}$$
and finally conclude
 $$\dim_{\mathcal{H}}(\mathcal{A})\geq t+(d-k)-\frac{(d-s)(\sigma-t)}{d+\sigma-(k+1)(d-k)}$$
by letting $t'\rightarrow t, s'\rightarrow s$, as desired.

\subsection{Cartesian products}\label{subsec-dual-Furstenberg-product}
Now we turn to the Cartesian product case  $E=E_1\times E_2\subset\R^{d-k}\times\R^k$. If we follow the argument above to reduce $E$ to a subset and then introduce a Frostman measure, the obstacle is to show this measure still has Cartesian product structure. This is why we need the extra assumption
$$0<\mathcal{H}^{s_1}(E_1), \mathcal{H}^{s_2}(E_2)<\infty.$$ The key fact here is, for every set $E$ there exists a Borel set $\tilde{E}$ 
satisfying 
\begin{equation}\label{Borel-closure}\mathcal{H}^{s}(\tilde{E})=\mathcal{H}^{s}(E),\ E\subset \tilde{E}\subset \mathcal{N}_\delta(E),\,\forall\,\delta>0,\end{equation} and if $\mathcal{H}^{s}(E)<\infty$ this $\tilde{E}$ is ``unique" in the sense that $\mathcal{H}^s(\tilde{E}\Delta\tilde{E}')=0$ for every other $\tilde{E}'$ satisfying \eqref{Borel-closure}. As this fact follows directly from the definition of Hausdorff measure, we leave details to the reader.

Now we construct Borel sets $\tilde{E_1}\subset\R^{d-k}, \tilde{E_2}\subset\R^k$ in this way and consider Frostman measures $\mu_1, \mu_2$ on them of exponents $s_1, s_2$ respectively. Then we run pigeonholing arguments similar to the above, to obtain $E''\subset E_1\times E_2$ such that $\mu_1\times\mu_2(\widetilde{E''})>0$ and
$$\mathcal{H}^{t'}_{\infty}(\{W+u\in\mathcal{A}:x\in W+u\})>c>0,\ \forall\,x\in E''.$$
Next we define $\gamma$, $j_0$, $\{B_i\}_{i\in I}$, $\mathcal{B}_j$, $\mathcal{B}_{j,x}$, $j_x$ in the same way and again write
$$E''=\bigcup_{j\geq j_0}\left\{x\in E'':j_x=j\right\}.$$
Now there is a small issue that it is not clear whether
\begin{equation}\label{small-issue}\widetilde{E''}=\bigcup_{j\geq j_0}\left(\widetilde{E''}\cap \widetilde{\left\{x\in E'':j_x=j\right\}}\right).\end{equation}
We point out that, although the two sides of \eqref{small-issue} may not equal as sets, they must have the same $\mu_1\times\mu_2$ measure. This is because, by the ``uniqueness" mentioned after \eqref{Borel-closure}, the symmetric difference between the two sides of \eqref{small-issue} must have $\mathcal{H}^{s_1}\times\mathcal{H}^{s_2}$ measure $0$, thus must be a $\mu_1\times\mu_2$ null set by the ball condition of Frostman measures. Therefore the pigeonhoing still works and there exists $j_1\geq j_0$ such that
\begin{equation}\label{j-1-square-product}\mu_1\times\mu_2(\widetilde{\{x\in E'':j_x=j_1\}})\gtrsim j_1^{-2}.\end{equation}
Then we can take
$$E_\delta:=\mathcal{N}_{\delta}(\{x\in E'':j_x=j_1\})$$
in the same way as above, and by ,\eqref{Borel-closure}\eqref{j-1-square-product} we have $$\mu_1\times\mu_2(E_\delta)\gtrsim j_1^{-2}.$$
So the same lower bound on the incidence \eqref{lower bound of incidence} still holds with $\mu_1\times\mu_2$. For the upper bound we just employ the second part of Theorem \ref{thm-incidence} in the same way as above with
$$p=2+\frac{s+\sigma-(k+1)(d-k)}{d-k-s_1}\left(\implies p'=1+\frac{d-k-s_1}{d-k+s_2+\sigma-(k+1)(d-k)}\right),$$
and eventually one can conclude that
$$\gamma\geq (t'-\sigma)p'+(d-k+\sigma)=t'+(d-k)-\frac{(d-k-s_1)(\sigma-t')}{d-k+s_2+\sigma-(k+1)(d-k)} $$
which implies
$$\dim_{\mathcal{H}}(\mathcal{A})\geq t+(d-k)-\frac{(d-k-s_1)(\sigma-t)}{d-k+s_2+\sigma-(k+1)(d-k)}$$
as desired. We leave details to readers as it is quite similar to the previous subsection.

\section{From incidence to discretized sum-product}\label{sec-sum-product}
In this section we prove Theorem \ref{thm-A+B-AC}. For convenience we denote by $|\cdot|$ the Lebesgue measure on $\R^2$, $\lambda_{2,1}$ on $\A(2,1)$, or their product.

For Theorem \ref{thm-A+B-AC}, we take $E=(A+B)\times(A\cdot C)$ and $$\mathcal{L}=\{x_2=c\cdot(x_1-b),\; b\in B, \; c\in C\}.$$ Then	$|A+B|\cdot|A\cdot C|\approx |E|$. 

For convenience we thicken $E$ to its $\delta$-neighborhood $F:=\mathcal{N}_{\delta}(E)$. As both $E, F$ are unions of $\delta$-balls, one can conclude $|E|\approx |F|$. Then we consider the incidence
	\begin{align*}
		I_{\delta}(F, \mathcal{L}):=\left\{(x,l)\in F\times\mathcal{L}:x\in\mathcal{N}_{\delta}(l)\right\}.
	\end{align*}

    First notice that, for every $l=\left\{(x_1,x_2):x_2=c\cdot(x_1-b)\right\}\in\mathcal{L}$,
    \begin{align*}
    	(a+b,a\cdot c)\in l\cap E,\ \forall\,a\in A.
    \end{align*}
    This implies
    \begin{align*}
    	|\left\{x\in F: x\in\mathcal{N}_{\delta}(l)\right\}|\gtrsim \delta\cdot|A|.
    \end{align*}
    As a result, we obtain a lower bound of the incidence
    \begin{align}\label{A+B-AC lower bound 1}
    	|I_{\delta}(F, \mathcal{L})|\gtrsim 
    	\delta\cdot |A|\cdot|\mathcal{L}|.
    \end{align}

    Then we turn to the upper bound. We would like to employ our estimate in the second half of Theorem \ref{thm-incidence}. To do so we need consider the  point-line duality as the non-concentration condition is on $B\times C$ but not necessarily on $(A+B)\times (A\cdot C)$.
    
     Take $\mathbf{D}: \mathbb{R}^2\to \mathbb{A}(2,1)$ as
	\begin{align*}
		(a,b)\mapsto x_2=ax_1+b:=l_{a,b},
	\end{align*}
    whose image contains all non-vertical lines. For the other direction consider the map $\mathbf{\tilde{D}}$ given by 
    \begin{align*}
    	l_{a,b}\mapsto (-a,b).
    \end{align*}
    Then $x\in l$ if and only if $\mathbf{D}(x)\in \mathbf{\tilde{D}}(l)$. As we only work on compact sets, both $\mathbf{D}, \mathbf{\tilde{D}}$ are bi-Lipschitz. Therefore it is equivalent to work on
    \begin{equation}\label{A+B-AC upper bound 1}
    |I_{c\delta}(\mathbf{\tilde{D}}(\mathcal{L}), \mathbf{D}(F)))|,
    \end{equation}
    where $c>0$ is a fixed constant. 

    To employ our estimates from Theorem \ref{thm-incidence} we need to check the non-concentration condition on $\mathbf{\tilde{D}}(\mathcal{L})$. By definition, $$\mathbf{\tilde{D}}(\mathcal{L})=\{(-c,-bc):b\in B, c\in C\},$$
    a bi-Lipschitz image of $B\times C$, thus a  $(\delta,s_B+s_C)$-set. Moreover, since the map $(x_1, x_2)\mapsto(-x_2, -x_1x_2)$ is bi-Lipschitz in $B\times C\subset[1,2]^2$, for every non-negative continuous function $f$,
    \begin{align}\label{desired-mu1-mu2-B-first}\int_{\mathbf{\tilde{D}}(\mathcal{L})}f(x)\,\frac{dx}{|\mathbf{\tilde{D}}(\mathcal{L})|}\approx & \int_C\int_B f(-x_2, -x_1x_2)\,\frac{dx_1}{|B|}\,\frac{dx_2}{|C|}\\\label{desired-mu1-mu2-C-first}=&\int_B\int_C f(-x_1, -x_1x_2)\,\frac{dx_1}{|C|}\,\frac{dx_2}{|B|}.\end{align}
    As a consequence, 
    $$\frac{1}{|\mathbf{\tilde{D}}(\mathcal{L})|}|I_{c\delta}(\mathbf{\tilde{D}}(\mathcal{L}), \mathbf{D}(F))|\approx \mu\times \lambda_{2,1}(I_{c\delta}(\mathbf{\tilde{D}}(\mathcal{L}), \mathbf{D}(F))),$$
    where $\mu$ is defined either as the right hand side of \eqref{desired-mu1-mu2-B-first} or \eqref{desired-mu1-mu2-C-first}. Notice that both satisfy \eqref{mu-product-like}. Say $\eqref{desired-mu1-mu2-B-first}$ for now.

    Now we are ready to invoke the second part of Theorem \ref{thm-incidence} with $d=2, k=\sigma=1$ to obtain
    \begin{equation}\label{invoke-thm-1.2}\frac{1}{|\mathbf{\tilde{D}}(\mathcal{L})|}|I_{c\delta}(\mathbf{\tilde{D}}(\mathcal{L}), \mathbf{D}(F))|\lesssim 
    I_s(\mu)^\frac{1}{p}\cdot \sup_{x_2} A_\alpha(\mu^{x_2}_1)^{\frac{p-2}{p}}\cdot\lambda_{2,1}(\mathcal{N}_{\delta}(\mathbf{D}(F)))^\frac{1}{p'}\cdot\delta, \end{equation}
    where $p:=2+\frac{s-1}{1-\alpha}$.

    By our discussion above, if we take $s=s_B+s_C-\epsilon$ and $\alpha=s_B-\epsilon$, then \eqref{invoke-thm-1.2} is reduced to
    $$\begin{aligned}|I_{c\delta}(\mathbf{\tilde{D}}(\mathcal{L}), \mathbf{D}(F))|\lesssim &  |\mathbf{\tilde{D}}(\mathcal{L})|\cdot\lambda_{2,1}(\mathcal{N}_{\delta}(D(F)))^\frac{1}{p'}\cdot\delta\\\approx &|\mathcal{L}|\cdot\lambda_{2,1}(\mathcal{N}_{\delta}(D(F)))^\frac{1}{p'}\cdot\delta\end{aligned},$$
    where $p'=\frac{1-s_B+s_C}{s_C}+O(\epsilon)$. Together with the lower bound \eqref{A+B-AC lower bound 1}, we conclude that
    $$|A+B|\cdot|A\cdot C|\approx \lambda_{2,1}(\mathcal{N}_{\delta}(D(F)))\gtrsim_{\epsilon, s_B, s_C} |A|^{p'}=|A|^{\frac{1-s_B+s_C}{s_C}+\epsilon}.$$

Finally, if we take $\mu$ as \eqref{desired-mu1-mu2-C-first}, the same argument leads to
$$|A+B|\cdot|A\cdot C|\gtrsim_{\epsilon, s_B, s_C} |A|^{\frac{1-s_C+s_B}{s_B}+\epsilon},$$
which completes the proof.

\bibliographystyle{abbrv}
\bibliography{mybibtex.bib}

\begin{thebibliography}{10}

\bibitem{Bou91}
J.~Bourgain.
\newblock Besicovitch type maximal operators and applications to {F}ourier
  analysis.
\newblock {\em Geom. Funct. Anal.}, 1(2):147--187, 1991.

\bibitem{Bou03}
J.~Bourgain.
\newblock On the {E}rd{\H o}s-{V}olkmann and {K}atz-{T}ao ring conjectures.
\newblock {\em Geom. Funct. Anal.}, 13(2):334--365, 2003.

\bibitem{DOV22}
D.~Dabrowski, T.~Orponen, and M.~Villa.
\newblock Integrability of orthogonal projections, and applications to
  {F}urstenberg sets.
\newblock {\em Adv. Math.}, 407:Paper No. 108567, 34, 2022.

\bibitem{Fal82}
K.~J. Falconer.
\newblock Hausdorff dimension and the exceptional set of projections.
\newblock {\em Mathematika}, 29(1):109--115, 1982.

\bibitem{FR24}
Y.~Fu and K.~Ren.
\newblock Incidence estimates for {$\alpha$}-dimensional tubes and
  {$\beta$}-dimensional balls in {$\Bbb R^2$}.
\newblock {\em J. Fractal Geom.}, 11(1-2):1--30, 2024.

\bibitem{Gan24}
S.~Gan.
\newblock Exceptional set estimate through {B}rascamp-{L}ieb inequality.
\newblock {\em Int. Math. Res. Not. IMRN}, (9):7944--7971, 2024.

\bibitem{GS77}
I.~M. Gel'fand and G.~E. Shilov.
\newblock {\em Generalized functions. {V}ol. 1}.
\newblock Academic Press [Harcourt Brace Jovanovich, Publishers], New
  York-London, 1964 [1977].
\newblock Properties and operations, Translated from the Russian by Eugene
  Saletan.

\bibitem{Hera19}
K.~H\'{e}ra.
\newblock Hausdorff dimension of {F}urstenberg-type sets associated to families
  of affine subspaces.
\newblock {\em Ann. Acad. Sci. Fenn. Math.}, 44(2):903--923, 2019.

\bibitem{HKM19}
K.~H\'{e}ra, T.~Keleti, and A.~M\'{a}th\'{e}.
\newblock Hausdorff dimension of unions of affine subspaces and of
  {F}urstenberg-type sets.
\newblock {\em J. Fractal Geom.}, 6(3):263--284, 2019.

\bibitem{KT01}
N.~H. Katz and T.~Tao.
\newblock Some connections between {F}alconer's distance set conjecture and
  sets of {F}urstenburg type.
\newblock {\em New York J. Math.}, 7:149--187, 2001.

\bibitem{Kau68}
R.~Kaufman.
\newblock On {H}ausdorff dimension of projections.
\newblock {\em Mathematika}, 15:153--155, 1968.

\bibitem{KM75}
R.~Kaufman and P.~Mattila.
\newblock Hausdorff dimension and exceptional sets of linear transformations.
\newblock {\em Ann. Acad. Sci. Fenn. Ser. A I Math.}, 1(2):387--392, 1975.

\bibitem{LWY21}
L.-H. Lim, K.~S.-W. Wong, and K.~Ye.
\newblock The {G}rassmannian of affine subspaces.
\newblock {\em Found. Comput. Math.}, 21(2):537--574, 2021.

\bibitem{Liu24}
B.~Liu.
\newblock Mixed-norm of orthogonal projections and analytic interpolation on
  dimensions of measures.
\newblock {\em Rev. Mat. Iberoam.}, 40(3):827--858, 2024.

\bibitem{Mar54}
J.~M. Marstrand.
\newblock Some fundamental geometrical properties of plane sets of fractional
  dimensions.
\newblock {\em Proc. London Math. Soc. (3)}, 4:257--302, 1954.

\bibitem{Mat95}
P.~Mattila.
\newblock {\em Geometry of Sets and Measures in Euclidean Spaces: Fractals and
  Rectifiability}, volume~44 of {\em Cambridge Studies in Advanced
  Mathematics}.
\newblock Cambridge University Press, Cambridge, 1995.

\bibitem{Mat15}
P.~Mattila.
\newblock {\em Fourier analysis and Hausdorff dimension}, volume 150.
\newblock Cambridge University Press, 2015.

\bibitem{PS00}
Y.~Peres and W.~Schlag.
\newblock Smoothness of projections, {B}ernoulli convolutions, and the
  dimension of exceptions.
\newblock {\em Duke Math. J.}, 102(2):193--251, 2000.

\bibitem{RW23}
K.~Ren and H.~Wang.
\newblock Furstenberg sets estimate in the plane.
\newblock {\em arXiv:2308.08819}, 2023.

\bibitem{Wolff99Kakeya}
T.~Wolff.
\newblock Recent work connected with the {K}akeya problem.
\newblock In {\em Prospects in mathematics ({P}rinceton, {NJ}, 1996)}, pages
  129--162. Amer. Math. Soc., Providence, RI, 1999.

\end{thebibliography}
\end{document}